# Phase retrieval from noisy data based on sparse approximation of object phase and amplitude

Vladimir Katkovnik

Signal Processing Laboratory, Technology University of Tampere, Tampere, Finland.

E-mail: vladimir.katkovnik@tut.fi.

*Abstract*—A variational approach to reconstruction of phase and amplitude of a complex-valued object from Poissonian intensity observations is developed. The observation model corresponds to the typical optical setups with a phase modulation of wavefronts. The transform domain sparsity is applied for the amplitude and phase modeling. It is demonstrated that this modeling results in the essential advantage of the developed algorithm for heavily noisy observations corresponding to a short exposure time in optical experiments. We consider also two simplified versions of this algorithm where the sparsity modeling of phase and amplitude is omitted. In the simulation study we compare the developed algorithms versus the Gerchberg-Saxton and truncation Wirtinger flow algorithms. The latter algorithm being the maximum likelihood based is the state-of-the-art for the phase retrieval from Poissonian observations. For noisy and very noisy observations the proposed algorithm demonstrates a valuable advantage.

## I. INTRODUCTION

### A. Phase retrieval formulation

ON many occasions a structure of specimens, for instance biological cells, is nearly invisible in intensity images. However, variations in thickness, density and refractive index result in variations of the phase delay. *Visualization* of these invisible phase variations by transforming them in light intensity is one of the challenging problems in optics. The revolutionary *phase contrast imaging* (Frits Zernike 1930s, Nobel prize 1953) solves the problem by introducing a modulation of the wavefront in the focal (Fourier) plane of the principal lens. In this way *visualization* is achieved as the observed light intensities are linked with the variations of phase. However, this is only a qualitative visualization because there is no one-to-one relations between the observed intensities and phase object properties. Nevertheless, the phase contrast microscopy is one of the frequently applied optical techniques in research and applications.

*Quantitative phase visualization* is targeted on precise phase imaging. It is fundamentally based on computational data processing. *Phase retrieval* is one of the computational techniques applied for quantitative phase imaging.

The $2D$ imaging is studied in this paper. Vectors and matrices in the following equations correspond to vectorized representations of 2D image variables. This vectorization is conventional in this kind of the problems.

The phase retrieval from multiple experiments is formulated as finding a complex-valued vector $\mathbf{x} \in \mathbb{C}^n$ from a set of real-valued observations $\mathbf{y}_s \in \mathbb{R}^m$:

$$\mathbf{y}_s = |\mathbf{A_s x}|^2, \ s = 1, ..., S. \qquad (1)$$

In terms of the coherent diffractive imaging the model (1) allows the following interpretation. Provided the unit intensity of a laser beam the vector $\mathbf{x}$ is an object (specimen) transfer function and in the same time it is a complex-valued wavefront just behind the object (see Fig.1a); $\mathbf{A}_s \in \mathbb{C}^{m \times n}$ is an $m \times n$ matrix of the wavefront propagation from the object to the sensor plane and the vector $\mathbf{y}_s$ is an intensity of the wavefront $\mathbf{u}_s = \mathbf{A_s x}$ registered by the sensor. The squared absolute value in (1) is an element-wise operation. Thus, the items of the vector $\mathbf{y}_s$ are squared absolute values of the corresponding items of the vector $\mathbf{A}_s \mathbf{x} \in \mathbb{C}^m$.

Assume for a moment that the complex-valued $\mathbf{u}_s = \mathbf{A}_s \mathbf{x}$ are known then the quadratic equations (1) are replaced by the linear ones and finding of $\mathbf{x}$ is reduced to the linear algebra problem. Thus, the phases of $\mathbf{u}_s$ eliminated by the modulus in (1) transform this linear problem in the quadratic problem which in general is much more complex.

Conventionally, the term *phase retrieval* is addressed to reconstruction of the missing phase in the vectors $\mathbf{u}_s$. However, the phase of the complex-valued object $\mathbf{x} \in \mathbb{C}^n$ is also unknown and actually reconstruction (retrieval) of this object phase is the main problem at hand.

In this paper we refocus the standard setting of the phase retrieval by treating the missed phases of $\mathbf{u}_s$ as auxiliary variables and the phase and the amplitude of the object $\mathbf{x}$ as the main variables of interest. We show that the sparse modeling applied to $\mathbf{x}$ (to both phase and amplitude) leads to an efficient algorithm solving the problem even for heavily noisy data.

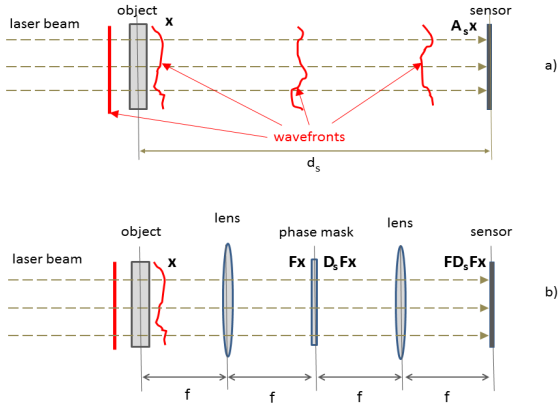

Fig. 1. Examples of optical setups for lensless (a) and $4f$ lens (b) scenarios for phase retrieval.

In this interpretation the phase retrieval becomes the quantitative phase imaging with the real-valued observations $\mathbf{y}_s \in \mathbb{R}^m$ given by (1).

Design of the image formation operators $\mathbf{A}_s$ in (1) is a crucial moment in order to gain an observation diversity sufficient for reliable reconstruction of the object phase as well as the object amplitude. The diversity means that the set $\{\mathbf{A}_s\}_1^S$ consists of the operators so different that the observations $\mathbf{y}_s$ suffice a finding $\mathbf{x}$. An image defocusing is one of the popular ways to get this kind of diversity. First results based on the defocus approach are demonstrated in [1] and [2] for two defocussed images obtain using an optical lens. A joint estimation of the object and the aberrations from multiple images incorporating the defocusing phase diversity is studied in [7] for an incoherent optical system.

In the coherent lensless imaging (Fig.1a) the laser beam goes through the object and after a free-space propagation of distance $d_s$ the intensity of the diffracted wavefield is registered by the sensor array. The corresponding operators $\mathbf{A}_s$ depend on the object-sensor distance $d_s$ and the wavelength $\lambda_s$, $\mathbf{A}_s = \mathcal{A}(d_s,\lambda_s)$, where the operator $\mathcal{A}$ can be modelled by the rigorous Rayleigh-Sommerfeld integral as well as by its Fresnel or Fraunhofer approximations [3]. These approximations can be calculated using the Fourier transform because their phase factors are cancelled by the modulus in (1). Then $\mathbf{A}_s = \mathbf{F}$, where $\mathbf{F}$ stays for the Fourier transform.

The approach using experiments with a set of distances $d_s$ (displaced sensor planes) is studied in [4], [5] and [6]. In the recent development, in particular we refer to [8] and [9], a spatial light modulator (SLM) is exploited in order to get a set of differently defocussed images.

A phase modulation of the wavefront near the object plane is another popular tool to gain the phase diversity. The phase modulation at the object plane plus the Fourier transform wavefront propagation result in the observation model known as a coded diffraction pattern [10]:

$$\mathbf{y}_s = |\mathbf{F}\mathbf{D}_s\mathbf{x}|^2, \ s=1,...,S, \quad (2)$$

where $\mathbf{F} \in \mathbb{C}^{n \times n}$ denotes the Fourier transform and $\mathbf{D}_s \in \mathbb{C}^{n \times n}$ is a diagonal matrix of complex exponents, $\mathbf{D}_s = diag\{\exp(j\phi_1(s)), \exp(j\phi_2(s)),...,\exp(j\phi_n(s))\}$.

This phase modulation can be implemented by a special phase mask inserted just behind the object plane in Fig.1a.

The phases $\phi_k(s)$ in $\mathbf{D}_s$ can be generated as random. Let $\phi_k(s)$ be i.i.d. zero mean Gaussian, $\phi_k(s) \sim N(0,\sigma)$. Then, $d_{ks} = \exp(j\phi_k(s))$ are random elements of $\mathbf{D}_s$, such that $E\{d_{ks}\} = 0$, here $E\{\cdot\}$ stays for the mathematical expectation, and the matrices $\mathbf{A}_s = \mathbf{F}\mathbf{D}_s$ are zero-mean random. This random phase modulation changes the object spectrum $\mathbf{F}\mathbf{x}$ in a radical way extending the spectrum intensity from low to higher frequencies.

The phase modulation (coded aperture imaging) is applied in various optical setups. For instance, let us consider the coherent $4f$ optical system shown in Fig.1b, where $f$ is a focal length of the lenses. The wavefront in the focal plane behind the first lens is the Fourier transform of the wavefront going through the object [3]. Further, the wavefront at the sensor plane, which locates at the focal length of the second lens, is the Fourier transform of the wavefront in the Fourier plane of the first lens. The phase modulation is produced using a phase mask or SLM in the focal (Fourier) plane of the first lens. The intensity registered by the sensor in this system can be represented as

$$\mathbf{y}_s = |\mathbf{F}\mathbf{D}_s\mathbf{F}\mathbf{x}|^2, \ s=1,...,S. \quad (3)$$

It is demonstrated in [11] that the phase $\phi_k(s)$ in $\mathbf{D}_s$ can be selected in such way that the wavefront at the sensor plane imitates desired displacements (defocus) $d_s$ of the sensor plane with respect to the object plane.

The setup shown in Fig.1a is a subject of an experimental study in this paper.

### B. Phase retrieval algorithms

Let us start from the popular Gerchberg-Saxton (GS) techniques (e.g. [12], [13]). These iterative algorithms are based on alternating projections between the object plane with a complex-valued $\mathbf{x}$ and the diffraction (Fourier) plane $\mathbf{Ax}$ with a given (measured) amplitude $\mathbf{z}$. At the diffraction plane the amplitudes in the vectors $\mathbf{Ax}$ are replaced by square roots of the corresponding items of $\mathbf{z}$. The back projection of this result to the object plane is modified according to the prior information on the object, e.g. support size and shape, amplitude value of $\mathbf{x}$, etc. The GS algorithms exist in many modifications. The review and analysis of the GS algorithms as well as further developments can be seen in the recent paper [14].

The algorithms known as a single-beam multiple-intensity reconstruction (SBMIR) are targeted on reconstruction of $3D$ wavefield covered by $2D$ intensity measurement planes (e.g. [4], [5], [6]). Citing these works: "The SBMIR algorithm starts from an initial guess of the wavefront at the first measurement plane. Then, this initial guess propagates numerically forward from the one measurement plane to the next following one successively through the all sequence of the measurements. At each plane the calculated modulus of the wave field is replaced by the square root of the intensity measured for this plane according to the GS algorithm. When the last measurement plane is reached the wave field estimate at this plane is propagated back to the first plane. This iterative process is repeated until convergence".

Contrary to the intuitively clear heuristic of the GS algorithms, the variational approaches to phase retrieval usually have a stronger mathematical background including image formation modeling, formulation of the objective function (criterion) and finally going to numerical techniques solving corresponding optimization tasks.

In particular, in incoherent imaging, the abberation-phase retrieval based on the maximum likelihood technique for Gaussian and Poissonian observations is developed in [7] with the phase diversity enabled by the multiple defocusing.

Here we wish refer to the recent overview [15] concentrated on the algorithms for the phase retrieval with the Fourier transform measurements of the form $\mathbf{y}=|\mathbf{Fx}|^2$ and some optical applications of these algorithms. The constrains sufficient for uniqueness of the solution are presented in detail. Beyond the alternating projection GS a few novel mathematical methods are discussed: semidefinite programming phase lifting using matrix completion *(PhaseLift algorithm)* [16] and greedy sparse phase retrieval *(GESPAR algorithm)* [17]. The fundamental progress for the methods based on the convex matrix optimization is announced in [18], where the novel algorithm is presented named "Sketchy Decisions". This algorithm allows to deal with high dimensions typical for the phase-lifting methods and is supported by the mathematical analysis with the convergence proof.

A sparse dictionary learning for phase retrieval is studied in [19]. A general phase retrieval algorithm to deal with noisy undersampled data corrupted by outliers is developed in [20]. The algorithm is based on a number of optimization problems solved with multiple initializations.

Many publications concern revisions of the intuitive GS algorithms by using optimization formulations. In particular, the links between the conventional GS and variational techniques are studied in [21] and [22]. A variational formulation for phase retrieval is demonstrated in [23], where the criterion for Poissonian observations and the prior imposing the phase smoothness are proposed. The problem is formalized as a penalized likelihood optimization. The conjugate gradient iterative algorithm for this setting is developed in [24].

In this brief overview, especially we wish to note the recent Wirtinger flow (WF) algorithms presented in [25] and [26]. These algorithms are iterative complex domain gradient descents applied to Poissonian likelihood criteria. Specific features of these algorithms are as follows: a spectral initialization, a non-trivial step-size parameter and the truncation of the gradient in the truncation Wirtinger flow (TWF) version of the algorithms [25]. The mathematical analysis is produced for the algorithm design, parameter selection and performance evaluation. It is stated that the solution of the quadratic equations (1) can be done "*nearly as easy as solving linear equation*". In this mathematical analysis the elements of the matrices $\mathbf{A}_s$ in (1) are random independent and subject of a complex-valued Gaussian distribution. Simulation experiments in [25] and [26] demonstrate that the TWF algorithm works and works very well for noiseless and low noise level observations.

The techniques based on the proximity operators developed in [20] and [27] provides a regularized optimization of the maximum likelihood criteria for Gaussian and Poissonian observations. The Bayesian approach to phase retrieval is developed in [28] for the Gaussian additive noise in the intensity observations.

The sparsity based techniques is a hot topic in phase retrieval. Publications based on the signal domain sparsity support minimization of the length of the vector-solution $\mathbf{x}$ (e.g. [17]). The transform domain amplitude and phase sparsity for the object $\mathbf{x}$ is developed for high-accuracy phase imaging (see [29], [30], [31]). In our paper [32] this sparsity modeling is applied for the phase retrieval in the 4f optical setup shown in Fig.1b.

The same sparsity modeling is exploited in this paper.

It is the only common point between this paper and [32] as they are different in optical setups, phase modulation methods, image formation and noise models and what is more important by a subject of research. This paper is oriented on phase retrieval from noisy and very noise Poissonian observations (photon-limited imaging) while in [32] the Gaussian observations are considered.

## C. Contribution and structure of this paper

In this paper we apply the adaptive nonlocal transform domain sparsity [29] for the phase retrieval from noisy coded diffraction patterns. The variational formulation of the algorithm design is produced for Poissonian observations. The derived Sparse Phase Amplitude Reconstruction (SPAR) algorithm has a structure typical for the GS style algorithms with forward and backward propagation of wavefronts. This algorithm incorporates two types of filtering: filtering of Poissonian observations at the sensor plane and filtering of phase and amplitude at the object plane. If these filterings are omitted the SPAR algorithm becomes similar to the conventional GS algorithm. Hereafter we use the term *GS algorithm* for this *simplified version* of the SPAR algorithm developed for phase retrieval.

Surprisingly, the GS algorithm demonstrates performance nearly identical to performance of the advanced TWF algorithm [26]. Both algorithms enable similar accuracy for phase and amplitude reconstruction as well as a similar computational complexity.

The SPAR algorithm computationally more demanding than GS demonstrates a much higher accuracy for noisy data as compared with respect to both TWF and GS. The phase unwrapping is included in the iterations of the SPAR algorithm when the object phase variation overcomes $2\pi$ range. It allows to achieve a more efficient noise suppression and a more accurate absolute phase reconstruction.

The random phase modulation of the wavefront and the transform domain sparsity modeling for phase and amplitude are essential features of the developed algorithm enabling the improved imaging for the high level noise in observations. The principal importance of the sparsity is demonstrated for undersampled observations where only the SPAR phase retrieval algorithm is efficient while other compared algorithms are failed.

The paper is organized as follows. In Section II the sparsity modeling for phase and amplitude and the Poissonian observations are discussed. The SPAR algorithm derivation is a subject of Section III, where step-by-step solutions of the variational problems are given for Poissonian observations. It is shown also that for the Gaussian noise in observations the respectively derived SPAR algorithm is different from the Poissonian SPAR algorithm only by filtering at the sensor plane. Section IV concerns the comparative experimental study of the algorithms for noisy Poissonian observations. The modification of the SPAR algorithm for the undersampled data and the corresponding experiments are presented in this section.

## II. PROBLEM FORMULATION

### A. Sparse wavefront modeling

It is recognized that many natural images (and signals) admit sparse representations, i.e. they can be well approximated by linear combinations of a small number of functions. This is a consequence of the so-called patch-wise self-similarity of these images. It means that is possible to find in them many quite similar small size patches located in different parts of the image. The sparse and redundant representations is of a special interest in the last years. This interest follows from importance that low dimensional models play in many signal and image applications such as compression, restoration, classification, just to mention a few of them [33].

Let $\mathbf{x} \in \mathbb{C}^n$ be a complex-valued wavefront. Denote $\mathbf{a} = \text{abs}(\mathbf{x})$ and $\boldsymbol{\varphi} = \text{angle}(\mathbf{x}) \in [-\pi, \pi)$ as, respectively, the corresponding images of amplitude (modulus) and the wrapped phase. Then we have $\mathbf{x} = \mathbf{a} \odot \exp(j\boldsymbol{\varphi})$, where '$\odot$' stands for the element-wise (Hadamard) product of two vectors. Herein, all functions applied to vectors are to be understood as component-wise.

With the intention of formulating treatable phase imaging problems, most approaches follow a two-step procedure: in the first step, an estimate of the so-called principal (wrapped, interferometric) phase in the interval $[-\pi, \pi)$ is determined; in the second step, termed phase unwrapping, the absolute phase is reconstructed by adding of an integer number of $2\pi$ multiples to the estimated interferometric phase [34]. In what follows, we denote the principal phase as $\boldsymbol{\varphi}$ and the absolute phase as $\boldsymbol{\varphi}_{abs}$. We introduce the phase-wrap operator $\mathcal{W}: \mathbb{R} \mapsto [-\pi, \pi)$, linking the absolute and principal phase as $\boldsymbol{\varphi} = \mathcal{W}(\boldsymbol{\varphi}_{abs})$. We also define the unwrapped phase as $\boldsymbol{\varphi}_{abs} = \mathcal{W}^{-1}(\boldsymbol{\varphi})$. Notice that $\mathcal{W}^{-1}$ is not an inverse operator for $\mathcal{W}$ because the latter is highly non-linear and for signals of dimension two and higher there is no one-to-one relation between $\boldsymbol{\varphi}_{abs}$ and $\boldsymbol{\varphi}$.

In sparse coding for complex-valued $\mathbf{x}$, we may think in two different directions: either we use a complex domain sparse representation to model directly for the complex image $\mathbf{x}$, as recently proposed in [37], [38] and [39], or we use separate sparse real-valued representations for the amplitude **a** and the absolute $\varphi_{abs}$ or interferometric $\varphi$ phase images of **x**.

The choice of the type of the sparse modeling depends on the application and a prior information. If the phase and amplitude are strongly correlated then the complex domain sparsity is preferable [39]. It should be noted that in the complex domain modeling the phase is treated only as a principal one.

In this paper, we follow the second type of the sparsity treating the phase and the amplitude as independent variables. Following to [29] we formalize this sparse wavefront modeling as the following matrix operations:

$$\mathbf{a} = \mathbf{\Psi}_a \boldsymbol{\theta}_a, \quad \boldsymbol{\varphi} = \mathbf{\Psi}_\varphi \boldsymbol{\theta}_\varphi, \tag{4}$$

$$\boldsymbol{\theta}_a = \mathbf{\Phi}_a \mathbf{a}, \quad \boldsymbol{\theta}_\varphi = \mathbf{\Phi}_\varphi \boldsymbol{\varphi}, \tag{5}$$

where $\boldsymbol{\theta}_a \in \mathbb{R}^p$ and $\boldsymbol{\theta}_\varphi \in \mathbb{R}^p$ are, respectively, amplitude and phase spectra of the object **x**. In (4), amplitude $\mathbf{a} \in \mathbb{R}^n$ and phase $\boldsymbol{\varphi} \in \mathbb{R}^n$ are synthesized from the amplitude and phase spectra $\boldsymbol{\theta}_a$ and $\boldsymbol{\theta}_\varphi$. On the other hand, the analysis Eqs.(5) give the spectra for amplitude and phase of the wavefront **x**. In (4)-(5) the synthesis $(n \times p)$ and analysis $(p \times n)$ matrices are denoted as $\mathbf{\Psi}_a, \mathbf{\Psi}_\varphi$ and $\mathbf{\Phi}_a, \mathbf{\Phi}_\varphi$, respectively.

Following the sparsity rationale we assume that the amplitude and phase spectra $\boldsymbol{\theta}_a$ and $\boldsymbol{\theta}_\varphi$ are sparse; i.e., most elements thereof are zero. In order to quantify the level of sparsity of $\boldsymbol{\theta}_a$ and $\boldsymbol{\theta}_\varphi$, i.e., their number of non-zero (active) elements, we use the pseudo $l_0$-norm $\|\cdot\|_0$ defined as a number of non-zero elements of the vector-argument. Therefore, we will design estimation criteria promoting low values of $\|\boldsymbol{\theta}_a\|_0$ and $\|\boldsymbol{\theta}_\varphi\|_0$.

Usually, the spectral dimensions are much higher than the dimensions of the image **x**, $p \gg n$, while the number of the active elements, i.e. the value of the pseudo $l_0$-norms of spectra, are much smaller than $p$ and sometimes can be even smaller than $n$.

It is obvious that for the complex exponent there is no difference between the principal and absolute phase, $\exp(j\varphi_{abs}) = \exp(j\varphi)$, and the angle operator in angle(**x**) gives the principal phase. However, there is a great deal of difference between the sparsity for absolute and interferometric phases. It is because in many applications the absolute phase can be smooth or piece-wise smooth function easily allowing sparsification while the corresponding wrapped phase experiences multiple discontinuities in particular if $\max(\mathrm{abs}(\varphi_{abs})) \gg \pi$. This kind of images are known as interferometric fringe patterns which are quite difficult for direct sparse modeling. In this case sparsification of phase through absolute phase modeling is preferable.

In what follows we treat the formulas (4)-(5) as universally applicable for principal and absolute phase. In the latter case $\varphi$ in (5) is replaced for $\varphi_{abs}$.

Note that in some cases an efficient sparsification of wrapped phase can be achieved through approximation of the complex exponent $\exp(j\varphi)$. Here we wish to mention the windowed Fourier transform developed for fringe processing in [40] and [41] as well as different forms of the Gabor transform which are definitely good candidates for this problem.

Another style of the data adaptive efficient approximators for the complex exponent are proposed in already mentioned papers [37], [38] and [39] based on the leaning dictionary techniques and high-order SVD non-local complex domain approximations.

### B. Noisy observation modeling

The measurement process in optics amounts to count the photons hitting the sensor's elements and is well modeled by independent Poisson random variables (e.g. [27], [42], [43]).

In many applications in biology and medicine the radiation (laser, X-ray, etc.) can be damaging for a specimen. Then, the dose (energy) of radiation can be restricted by a lower exposure time or by use a lower power radiation source, say up to a few or less numbers of photons per pixel of sensor what leads to heavily noisy registered measurements. Imaging from these observations, in particular, phase imaging is called photon-limited.

The probability that a random Poissonian variable $\mathbf{z}_s[l]$ of the mean value $\mathbf{y}_s[l]$ takes a given non-negative integer value $k$, is given by

$$p(\mathbf{z}_s[l] = k) = \exp(-\mathbf{y}_s[l]\chi)\frac{(\mathbf{y}_s[l]\chi)^k}{k!}, \tag{6}$$

where $\mathbf{y}_s[l]$ is the intensity of the wavefront at the pixel $l$ (1).

Recall that the mean and the variance of Poisson random variable $\mathbf{z}_s[l]$ are equal and are given by $\mathbf{y}_s[l]\chi$, i.e., $E\{\mathbf{z}_s[l]\} = \mathrm{var}\{\mathbf{z}_s[l]\} = \mathbf{y}_s[l]\chi$. Defining the observation signal-to-noise ratio (SNR) as the ratio between the square of the mean and the variance of $\mathbf{z}_s[l]$, we have $SNR = E^2\{\mathbf{z}_s[l]\}/\mathrm{var}\{\mathbf{z}_s[l]\} = \mathbf{y}_s[l]\chi$. Thus, the relative noisiness of observations becomes stronger as $\chi \to 0$ ($SNR \to 0$) and approaches zero when $\chi \to \infty$ ($SNR \to \infty$). The latter case corresponds to the noiseless scenario, with the probability 1

$$\mathbf{z}_s[l]/\chi \to \mathbf{y}_s[l]. \tag{7}$$

The parameter $\chi > 0$ in (6) is a scaling factor defining a proportion between the intensity of the observations with respect to the intensity of the input wavefront. This parameter is of importance as it controls a level of the noise in observations. Physically it can be interpreted as an exposure time and as the sensitivity of the sensor with respect to the input radiation.

In order to make the noise more understandable the noise level can be characterized by the estimates of SNR

$$SNR = 10\log_{10}(\chi^2 \sum_{s=1}^{S}||\mathbf{y}_s||_2^2 / \sum_{s=1}^{S}||\mathbf{y}_s\chi - \mathbf{z}_s||_2^2) \ dB \tag{8}$$

and of the mean value of photons per pixel:

$$N_{photon} = \sum_{s=1}^{S}\sum_{l=1}^{n}\mathbf{z}_s[l]/Sn. \tag{9}$$

A smaller values of $\chi$ lead to smaller $SNR$ and $N_{photon}$, i.e. to noisier observations $\mathbf{z}_s$.

## III. ALGORITHM DEVELOPMENT

We consider the problem of wavefront reconstruction as an estimation of $\mathbf{x} \in \mathbb{C}^n$ from the Poissonian observations $\{\mathbf{z}_s\}_1^S$. This problem is rather challenging mainly due the periodic nature of the likelihood function with respect to the phase $\varphi$ and the non-linearity of the observation model.

The maximum likelihood means minimization on $\mathbf{x}$ of the negative log-likelihood criterion:

$$\mathcal{L}(\{\mathbf{u}_s\}) = \sum_{s=1}^{S}\sum_{l=1}^{n}[|\mathbf{u}_s[l]|^2\chi - \mathbf{z}_s[l]\log(|\mathbf{u}_s[l]|^2\chi)] \tag{10}$$

corresponding to the observations (6).

The WF and TWF algorithms in [25] and [26] implement minimization of (10) based on calculation of the gradient of $\partial \mathcal{L}(\{\mathbf{A}_s\mathbf{x}\})/\partial \mathbf{x}$, $\mathbf{x} \in \mathbb{C}^n$, and on the corresponding gradient descent iterations.

Herein contrary to this straightforward approach, we adopt multiobjective Nash equilibrium optimization. The main objective of this approach is a simultaneous minimization of the negative log-likelihood function (10) and of the $l_0$-norms of the magnitude and phase spectra given as $||\boldsymbol{\theta}_a||_0$ and $||\boldsymbol{\theta}_\varphi||_0$, respectively. However, the approach based on these intentions yields complex calculations with respect to $(\mathbf{a}, \varphi)$.

In order to make the problem manageable, we introduce auxiliary variables $\mathbf{v}_s$ approximating the wavefronts $\mathbf{u}_s$ and allowing to split optimization with respect to $(\mathbf{a}, \varphi, \boldsymbol{\theta}_a, \boldsymbol{\theta}_\varphi)$ into simpler decoupled problems. These successive steps are introduced in the following subsection.

In the Nash equilibrium formulation, as it is implemented in this paper, the conventional constrained optimization with a single criterion function is replaced by balancing two criteria. Details of this approach, links with the game theory and demonstrations of its efficiency for the synthesis-analysis sparse inverse imaging can be seen in [36], where it is done for linear real-valued observations. The Nash equilibrium technique for complex domain problems in optics was applied in [29]-[32], [44].

### A. SPAR algorithm

The following two criteria are introduced for formalization of the algorithm design [44]:

$$\mathcal{L}_1(\{\mathbf{u}_s\}, \mathbf{x}) = \sum_{s=1}^{S}\sum_{l=1}^{n}[|\mathbf{u}_s[l]|^2\chi - \mathbf{z}_s[l]\log(|\mathbf{u}_s[l]|^2\chi)] + \frac{1}{\gamma_1}\sum_{s=1}^{S}||\mathbf{u}_s - \mathbf{A}_s\mathbf{x}||_2^2, \tag{11}$$

$$\mathcal{L}_2(\boldsymbol{\theta}_\varphi, \boldsymbol{\theta}_a, \mathbf{x}) = \tau_a \cdot ||\boldsymbol{\theta}_a||_0 + \tau_\varphi \cdot ||\boldsymbol{\theta}_\varphi||_0 + \frac{1}{2}||\boldsymbol{\theta}_a - \Phi_a\mathbf{a}||_2^2 + \frac{1}{2}||\boldsymbol{\theta}_\varphi - \Phi_\varphi\boldsymbol{\varphi}||_2^2. \tag{12}$$

The criterion (11) is a regularized version of (10) with the quadratic penalization by the differences between the wavefronts $\mathbf{u}_s$ and their estimates $\mathbf{v}_s = \mathbf{A}_s\mathbf{x}$ used as splitting variables.

It is noted in Subsection II-A that we use the separate sparse modeling for the phase $\varphi$ and the amplitude $\mathbf{a}$ of the wavefront $\mathbf{x}$. The criterion (12) promotes this sparsity in the transform domain. The regularization terms $\frac{1}{2}||\boldsymbol{\theta}_a - \Phi_a\mathbf{a}||_2^2$ and $\frac{1}{2}||\boldsymbol{\theta}_\varphi - \Phi_\varphi\boldsymbol{\varphi}||_2^2$ are squared Euclidean norms calculated for differences between spectra $\boldsymbol{\theta}_a$ and $\boldsymbol{\theta}_\varphi$ and their predictors $\Phi_a\mathbf{a}$ and $\Phi_\varphi\boldsymbol{\varphi}$.

Accordingly to the very idea of the Nash equilibrium for balancing multiple penalty functions (e.g. [45]) the proposed algorithm is composed of alternating optimization steps performed for the criteria $\mathcal{L}_1$-$\mathcal{L}_2$. It leads to the iterative algorithm:

$$\{\hat{\mathbf{u}}_s^t\} = \arg\min_{\{\mathbf{u}_s\}} \mathcal{L}_1(\{\mathbf{u}_s\}, \hat{\mathbf{x}}^t), \tag{13}$$

$$\hat{\mathbf{x}}^t = \arg\min_{\mathbf{x}} \mathcal{L}_1(\{\hat{\mathbf{u}}_s^t\}, \mathbf{x}), \tag{14}$$

$$(\hat{\boldsymbol{\theta}}_\varphi^t, \hat{\boldsymbol{\theta}}_a^t) = \arg\min_{\boldsymbol{\theta}_\varphi, \boldsymbol{\theta}_a} \mathcal{L}_2(\boldsymbol{\theta}_\varphi, \boldsymbol{\theta}_a, \hat{\mathbf{x}}^t), \tag{15}$$

$$\hat{\mathbf{a}}^{t+1} = \Psi_a \hat{\boldsymbol{\theta}}_a^t, \ \hat{\boldsymbol{\varphi}}^{t+1} = \Psi_\varphi \hat{\boldsymbol{\theta}}_\varphi^t \tag{16}$$

$$\hat{\mathbf{x}}^{t+1} = \hat{\mathbf{a}}^{t+1} \odot \exp(j\hat{\boldsymbol{\varphi}}^{t+1}), \tag{17}$$

where the last two equation (16)-(17) update amplitude, phase and complex-valued wavefront: $\hat{\mathbf{a}}^{t+1}$, $\hat{\boldsymbol{\varphi}}^{t+1}$, $\hat{\mathbf{x}}^{t+1}$.

Let us solve optimizations in (13)-(15).

(1) Note that $\mathcal{L}_1(\{\mathbf{u}_s\}, \mathbf{x})$ is additive with respect to $\mathbf{u}_s[l]$ and minimization on $\{\mathbf{u}_s\}$ (the problem (13))

is reduced to the scalar minimization on $\mathbf{u}_s[l]$ of the criterion

$$|\mathbf{u}_s[l]|^2\chi - \mathbf{z}_s[l]\log(|\mathbf{u}_s[l]|^2\chi) + \frac{1}{\gamma_1}|\mathbf{u}_s[l]-\mathbf{v}_s[l]|^2. \quad (18)$$

It can be verified that this minimization on the complex-valued $\mathbf{u}_s[l]$ gives angle($\mathbf{u}_s[l]$)=angle($\mathbf{v}_s[l]$) [44].

Inserting this solution in (18) we obtain the criterion depending only on the absolute values $|\mathbf{u}_s[l]|$:

$$|\mathbf{u}_s[l]|^2\chi - \mathbf{z}_s[l]\log(|\mathbf{u}_s[l]|^2\chi) + \frac{1}{\gamma_1}||\mathbf{u}_s[l]|-|\mathbf{v}_s[l]||^2. \quad (19)$$

After differentiation on $|\mathbf{u}_s[l]|$ we obtain the quadratic equation for the optimal absolute values. A unique non-negative root of this quadratic equation is [44]:

$$\mathbf{b}_s[l] = \frac{|\mathbf{v}_s[l]| + \sqrt{|\mathbf{v}_s[l]|^2 + 4\mathbf{z}_s[l]\gamma_1(1+\gamma_1\chi)}}{2(1+\gamma_1\chi)}. \quad (20)$$

Then, the solution for (13) is of the form:

$$\mathbf{u}_s[l] = \mathbf{b}_s[l]\exp(j \cdot \text{angle}(\mathbf{v}_s[l])). \quad (21)$$

In this solution the amplitude $\mathbf{b}_s[l]$ depends on both the observation $\mathbf{z}_s[l]$ and the amplitude of $\mathbf{v}_s[l]$.

It is useful to note that minimization of $\mathcal{L}_1(\{\mathbf{u}_s\},\mathbf{x})$ on $\{\mathbf{u}_s\}$ can be treated as the proximity operator for Poissonian observations as it is presented in [20] and [27].

For the large $\chi \to \infty$ (noiseless case)

$$\mathbf{u}_s[l] \to \sqrt{\mathbf{z}_s[l]/\chi}\exp(j \cdot angle(\mathbf{v}_s[l])). \quad (22)$$

According to (7), $\mathbf{z}_s[l]/\chi \to \mathbf{y}_s[l]$ and the amplitude update formula is as follows

$$\mathbf{u}_s[l] = \sqrt{\mathbf{y}_s[l]}\exp(j \cdot \text{angle}(\mathbf{v}_s[l])), \ s=1,...,S. \quad (23)$$

(2) Optimization of $\mathcal{L}_1(\{\mathbf{u}_s^t\},\mathbf{x})$ with respect to $\mathbf{x}\in\mathbb{C}^n$ (the problem (14)) leads to the minimum condition of the form $\partial\mathcal{L}_1(\{\mathbf{u}_s^t\},\mathbf{x})/\partial\mathbf{x}^* = 0$ and to the normal least-squares equation for $\mathbf{x}$

$$\sum_{s=1}^S \mathbf{A}_s^H\mathbf{A}_s\mathbf{x} = \sum_{s=1}^S \mathbf{A}_s^H\mathbf{u}_s \quad (24)$$

and to the solution

$$\mathbf{x} = (\sum_{s=1}^S \mathbf{A}_s^H\mathbf{A}_s)^{-1}\sum_{s=1}^S \mathbf{A}_s^H\mathbf{u}_s, \quad (25)$$

provided that the matrix $\sum_{s=1}^S \mathbf{A}_s^H\mathbf{A}_s$ is non-singular.

In particular, for the Fourier forward propagation $\mathbf{A}_s = \mathbf{D}_s\mathbf{F}$ and $\mathbf{A}_s^H\mathbf{A}_s = \mathbf{I}_{n\times n}$ provided $n=m$. Then (25) takes the form

$$\mathbf{x} = \frac{1}{S}\sum_{s=1}^S \mathbf{A}_s^H\mathbf{u}_s. \quad (26)$$

In general, the situation may be more complex, in particular, because $\mathbf{A}_s^H\mathbf{A}_s$ are ill-conditioned due to the fact that the operators $\mathbf{A}_s$ are low-pass filters suppressing high frequency components of the object $\mathbf{x}$ [3] or because of undersampled observations, $m < n$.

Then, the solution of (24) can be found using one of the methods applicable to the normal equation (24) (see e.g. [46]).

Note, that in (24) $\mathbf{A}_s\mathbf{e}^k$ stays for the forward propagation of the wavefront $\mathbf{e}^k$ and $\mathbf{A}_s^H\mathbf{A}_s\mathbf{e}^k$ means the backward propagation of the wavefront $\mathbf{A}_s\mathbf{e}^k$. It follows that the algorithm for solution of (24) can implemented bypassing the large size matrices $\mathbf{A}_s$ by use of forward and backward propagation operators applied to $2D$ images.

(3) In general, multidimensional minimization of the $l_0$-pseudonorm results in the hard computational problem. The considered criterion $\mathcal{L}_2(\boldsymbol{\theta}_\varphi, \boldsymbol{\theta}_a, \mathbf{x})$ (the problem (15)) is additive with respect to the items of the vectors $\boldsymbol{\theta}_a, \boldsymbol{\theta}_\varphi$. Then the minimization can be produced independently for each items of these vectors and this scalar optimization allows a simple analytical solution:

$$\hat{\boldsymbol{\theta}}_a = (\Phi_a\mathbf{a})\odot 1\left[\text{abs}(\Phi_a\mathbf{a}) \geq \sqrt{2\tau_a}\right], \quad (27)$$
$$\hat{\boldsymbol{\theta}}_\varphi = (\Phi_\varphi\boldsymbol{\varphi})\odot 1\left[\text{abs}(\Phi_\varphi\boldsymbol{\varphi}) \geq \sqrt{2\tau_\varphi}\right],$$

where $1[\mathbf{w}]$, $\mathbf{w}\in\mathbb{R}^p$, is an element-wise vector function: $\mathbb{R}^p \mapsto \mathbb{R}^p$, $1[\mathbf{w}_k] = 1$ if $\mathbf{w}_k \geq 0$ and $1[\mathbf{w}_k] = 0$ if $\mathbf{w}_k < 0$.

The formulas (27) define the thresholding (hard-thresholding) operation. Here $th_a = \sqrt{2\tau_a}$ and $th_\varphi = \sqrt{2\tau_\varphi}$ are the thresholding parameters for the amplitude and the phase, respectively. The items of the spectral coefficients abs($\Phi_a\mathbf{a}$) and abs($\Phi_\varphi\boldsymbol{\varphi}$), which are smaller than the corresponding thresholds are zeroed in (27).

The success of any sparse imaging depends on how reach and redundant are transforms/dictionaries used for analysis and synthesis. In our algorithm for the analysis and synthesis operations we use the BM3D frames (matrices), where BM3D is the abbreviation for Block-Matching and 3D filtering [35], [36].

Let us recall some basic ideas of this popular technique. At the first stage the image is partitioned into small overlapping square patches. For each patch a group of similar patches is collected which are stacked together and form a 3D array (group). This stage is called *grouping*. The entire 3D group-array is projected onto a 3D transform basis. The obtained spectral coefficients are hard-thresholded and the inverse 3D transform gives the filtered patches, which are returned to the original position of these patches in the image. This stage is called *collaborative* filtering. This process is repeated for all pixels of the entire wavefront and obtained overlapped filtered patches are aggregated in the final image estimate. This last stage is called *aggregation*. The details of BM3D as an advanced image filter can be seen in [35].

It follows from [36] and [44], that the steps (15)-(16) including the *grouping* operations defining the analysis $\Phi$ and synthesis $\Psi$ matrices can be combined in a single algorithm. In what follows, we use the notation BM3D for this algorithm. Note, that the standard BM3D algorithm as it is presented in the original paper [35] is composed from two successive stages: thresholding and Wiener filtering. The BM3D algorithm corresponding to the procedures (15) and (16) consists of the first thresholding stage only.

The criterion $\mathcal{L}_2$ is separable on $\boldsymbol{\theta}_\varphi$ and $\boldsymbol{\theta}_a$. It follows that the corresponding solutions can be calculated independently for amplitude and phase. Using the BM3D algorithm for implementation of the steps (15)-(16) we obtain:

$$\begin{aligned}\hat{\mathbf{a}}^{t+1} &= BM3D_{ampl}(\hat{\mathbf{a}}^t, th_a), \\ \hat{\boldsymbol{\varphi}}^{t+1} &= BM3D_{phase}(\hat{\boldsymbol{\varphi}}^t, th_\varphi).\end{aligned} \quad (28)$$

In (28) BM3D stays with different subscripts because different parameters can be used in BM3D for amplitude and phase processing.

Combining the solutions obtain for (13)-(14) and (28) as the solution for (15)-(16) we arrive to the phase retrieval algorithm shown in Table I.

TABLE I
SPAR PHASE RETRIEVAL ALGORITHM

|   | Input: $\{\mathbf{z}_s\}$, $s = 1, ..., S$, $\mathbf{x}^1$; |
|---|---|
|   | For $t = 1, ..., N$; |
| 1. | **Forward propagation**: |
|   | $\hat{\mathbf{v}}_s^t = \mathbf{A}_s \hat{\mathbf{x}}^t$, $s = 1, ..., S$; |
| 2. | **Poissonian noise filtering**: |
|   | $\hat{\mathbf{u}}_s^t = \hat{\mathbf{b}}_s^t \odot \exp(j \cdot angle(\hat{\mathbf{v}}_s^t))$, Eq.(20) for $\hat{\mathbf{b}}_s^t$; |
| 3. | **Backward propagation**: |
|   | $\hat{\mathbf{x}}^t = (\sum_{s=1}^S \mathbf{A}_s^H \mathbf{A}_s)^{-1} \sum_{s=1}^S \mathbf{A}_s^H \hat{\mathbf{u}}_s^t$; |
| 4. | **Phase unwrapping**: |
|   | $\hat{\boldsymbol{\varphi}}_{abs}^t = \mathcal{W}^{-1}(angle(\hat{\mathbf{x}}^t))$; |
| 5. | **Phase and amplitude filtering**: |
|   | $\hat{\boldsymbol{\varphi}}_{abs}^{t+1} = BM3D_{phase}(\hat{\boldsymbol{\varphi}}_{abs}^t, th_\varphi)$, |
|   | $\hat{\mathbf{a}}^{t+1} = BM3D_{ampl}(abs(\hat{\mathbf{x}}^t), th_a)$; |
| 6. | **Object wavefront update**: |
|   | $\hat{\mathbf{x}}^{t+1} = \mathbf{a}^{t+1} \circ \exp(j\hat{\boldsymbol{\varphi}}_{abs}^{t+1})$; |
|   | Output: $\hat{\boldsymbol{\varphi}}_{abs}^{N+1}$, $\hat{\mathbf{a}}^{N+1}$. |

At Step 1 the object wavefront estimate $\hat{\mathbf{x}}^t$ propagates using the operators $\mathbf{A}_s$ and defines the wavefront $\hat{\mathbf{v}}_s^t$ at the sensor plane. At Step 2 this wavefront is updated to the variable $\hat{\mathbf{u}}_s^t$ by changing the amplitude according to the given observations $\mathbf{z}_s$ (20), while the phase of $\hat{\mathbf{v}}_s^t$ is preserved in $\hat{\mathbf{u}}_s^t$. At Step 3 the estimates $\{\hat{\mathbf{u}}_s^t\}$ backpropagate to the object plane and update the object wavefront $\hat{\mathbf{x}}^{t+1}$.

The sparsification (filtering on the base of sparse approximations) is produced in Step 5. The algorithm in Table I is presented in the form for sparsification of the absolute phase. In order to get the estimate of the absolute phase the unwrapping operation is given as Step 4.

If the phase sparsity is imposed on the principal phase then Step 4 is omitted and $\hat{\boldsymbol{\varphi}}_{abs}^t = angle(\hat{\mathbf{x}}^t)$ in the following Step 5.

The derived SPAR can be modified for the Gaussian noise in observations, i.e. for

$$\mathbf{z}_s[l] = |\mathbf{u}_s[l]|^2 + \sigma \boldsymbol{\varepsilon}_s[l] \quad (29)$$

where $\varepsilon_s$ is i.i.d. Gaussian errors, $\varepsilon_s[l] \sim \mathcal{N}(0,1)$, and $\sigma$ stands for the standard deviation of the noise.

Then criterion $\mathcal{L}_1$ (11) takes the form [32]

$$\mathcal{L}_1(\{\mathbf{u}_s\}, \mathbf{x}) = \frac{1}{\sigma^2}\sum_{s=1}^S\sum_{l=1}^n[|\mathbf{u}_s[l]|^2 - \mathbf{z}_s[l]]^2 + \quad (30)$$

$$\frac{1}{\gamma_1}\sum_{s=1}^S ||\mathbf{u}_s - \mathbf{v}_s||_2^2.$$

The procedure (13)-(17) with $\mathcal{L}_1$ is defined by (30) gives the algorithm in the form identical to shown in Table I with the only difference in Step 2 where the amplitude update $\hat{\mathbf{b}}_s[l]$ is calculated as the solution of the Cardan equation [32] (see also [27]):

$$\mathbf{b}_s^3[l] + C\mathbf{b}_s[l] + D = 0, \quad (31)$$
$$C = \frac{\sigma^2}{2\gamma_1} - z_s[l], \ D = -\frac{\sigma^2}{2\gamma_1}|\mathbf{v}_s[l]|.$$

As the coefficient at the second degree in (31) is equal to zero and $D \leq 0$ this equation has a unique real non-negative solution which locates between $\sqrt{z_s[l]}$ and $|\mathbf{v}_s[l]|$. If $\gamma_1 \to \infty$ this root approaches $\sqrt{z_s[l]}$ and if $\gamma_1 \to 0$ this root approaches $|\mathbf{v}_s[l]|$.

### B. GS algorithms

As a particular case, we introduce a simplified version of the SPAR algorithm where Steps 4 and 5 are dropped (see Table II). In this case the filtering of phase and amplitude, i.e. sparsity modeling, are excluded while the noise filtering in Step 2 is preserved.

We name this algorithm GS-F, where "F" reminds that the observation filtering is applied, and use the name GS for the algorithm where the observation filtering is also dropped and Step 2 is implemented accordingly to the standard rule (22). In our derivation it corresponds to the noiseless data. This latter GS algorithm corresponds to the conventional GS being applied for the phase retrieval from coded diffraction patterns.

Note that contrary to the conventional GS heuristic, here the algorithms the GS-F and GS algorithms are derived from minimization of $\mathcal{L}_1(\{\mathbf{u}_s\}, \mathbf{x})$ with respect to $\{\mathbf{u}_s\}$ and $\mathbf{x}$. The parameter $\gamma_1$ in the updating rule (20) is responsible for the filtering properties of the GS-F algorithm: smaller $\gamma_1$ means the stronger filtering of the observations.

TABLE II
GS-F PHASE RETRIEVAL ALGORITHM

| | Input: $\{\mathbf{z}_s\}$, $s=1,...,S$, $\mathbf{x}^1$; |
|---|---|
| | For $t=1,...,N$; |
| 1. | **Forward propagation**: |
| | $\hat{\mathbf{v}}_s^t = \mathbf{A}_s \hat{\mathbf{x}}^t$, $s=1,...,S$; |
| 2. | **Observation constrains**: |
| | $\hat{\mathbf{u}}_s^t = \hat{\mathbf{b}}_s^t \odot \exp(j \cdot angle(\hat{\mathbf{v}}_s^t))$, Eq.(20) for $\hat{\mathbf{b}}_s^t$; |
| 3. | **Backward propagation**: |
| | $\hat{\mathbf{x}}^{t+1} = (\sum_{s=1}^S \mathbf{A}_s^H \mathbf{A}_s)^{-1} \sum_{s=1}^S \mathbf{A}_s^H \hat{\mathbf{u}}_s^t$; |
| | Output: $\hat{\mathbf{x}}^{N+1}$. |

*C. Algorithm's implementation*

In implementation of the SPAR algorithm the vectorization is bypassed and calculations are produced for $2D$ variables. It becomes possible mainly due to the corresponding implementation of BM3D where the analysis and syntheses operations are produced algorithmically without formation of the analysis and syntheses matrices $\Phi$ and $\Psi$. The implementation of the SPAR algorithm for the coded diffraction pattern scenario (2) is shown in Table III, where we preserve the notation of the variable but replace the bold fonts used for the vectorized representation by the corresponding normal ones.

Here $D_s$ is a phase modulation mask (image) composed from complex exponents used in (2). It should be emphasized that the analysis and synthesis transforms (frames) $\Phi$ and $\Psi$, as they are introduced in (4) and (5), are varying from iteration to iteration because the grouping in BM3D depends on the input image and on each iterations these images are different.

In our implementation the thresholds $th_\varphi$ and $th_a$ in BM3D filters are data adaptive. They are calculated as

$$th_\varphi = th_\varphi^0 \cdot \hat{\sigma}(\hat{\varphi}_{abs}^t), \ th_a = th_a^0 \cdot \hat{\sigma}(abs(\hat{x}^t)), \quad (32)$$

where $th_\varphi^0$ and $th_a^0$ are invariant parameters and $\hat{\sigma}(\hat{\varphi}_{abs}^t)$ and $\hat{\sigma}(abs(\hat{x}^t))$ are estimates of the standard deviation of the variables-estimates $\hat{\varphi}_{abs}^t$ and $abs(\hat{x}^t)$, respectively. Thus, the smoothing properties of the BM3D filters become stronger for larger standard deviations and softer for smaller values of standard deviations.

TABLE III
IMPLEMENTATION OF SPAR PHASE RETRIEVAL ALGORITHM

| | Input: $\{z_s\}$, $s=1,...,S$, $x^1$; |
|---|---|
| | For $t=1,...,N$; |
| 1. | **Forward propagation**: |
| | $\hat{v}_s^t = \mathcal{F}\{D_s x^t\}$, $s=1,...,S$; |
| 2. | **Poissonian noise suppression**: |
| | $\hat{u}_s^t = \hat{b}_s^t \exp(j \cdot angle(\hat{v}_s^t))$, Eq.(20) for $\hat{b}_s^t$; |
| 3. | **Backward propagation**: |
| | $\hat{x}^t = \frac{1}{S} \sum_{s=1}^S \mathcal{F}^{-1}\{D_s^* \hat{u}_s^t\}$; |
| 4. | **Phase unwrapping**: |
| | $\hat{\varphi}_{abs}^t = \mathcal{W}^{-1}(angle(\hat{x}^t))$; |
| 5. | **Phase and amplitude filtering**: |
| | $\hat{\varphi}_{abs}^{t+1} = BM3D_{phase}(\hat{\varphi}_{abs}^t, th_\varphi)$, |
| | $\hat{b}^{t+1} = BM3D_{ampl}(abs(\hat{x}^t), th_a)$; |
| 6. | **Object wavefront update**: |
| | $\hat{x}^{t+1} = \hat{a}^{t+1} \odot \exp(j\hat{\varphi}_{abs}^{t+1})$; |
| | Output: $\hat{\varphi}_{est} = \hat{\varphi}_{abs}^{N+1}$, $\hat{a}_{est} = \hat{a}^{N+1}$. |

For phase unwrapping we exploit the PUMA algorithm [47] based on an energy minimization by the graph cut techniques. The PUMA is able to minimize a wide class of energies, defining flexibility of the method. It is one of the best phase unwrapping algorithm in the field with a unique ability to reconstruct discontinuous absolute phases.

IV. NUMERICAL EXPERIMENTS

For simulation tests we select the coded diffraction pattern scenario (2) with the algorithm's implementation shown in Table III.

Following to the publications [25] and [26] the wavefront modulation is enabled by the random phases $\phi_k$ in $D_s$ with equal probabilities taking four values $[0, \pi/2, -\pi/2, \pi]$. The choice of these four random phase values for phase modulation is caused by our intention to consider TWF as a main counterpart to our algorithms. The MATLAB codes of TWF provided by the authors make the comparative analysis simple for implementation. While the experimental tests in [26] are presented mainly for noiseless data or small level noise herein we are concentrated on noisy data and show that in this case the sparse modeling developed in this paper allows to achieve a dramatic improvement in the accuracy of phase and amplitude imaging.

In the phase retrieval problem the object phase image can be estimated within an invariant phase-shift only. Following [25] and [26] the estimated principal phase image is corrected by an invariant phase-shift $\varphi_{shift}$ defined as

$$\varphi_{shift} = \arg(\min_{\varphi \in [0,2\pi]} ||\exp(-j\varphi)\hat{x}_{est} - x_{true}||_2^2), \quad (33)$$

here $x_{true}$ and $\hat{x}_{est}$ are the true wavefront and the estimate, respectively. This correction of the phase is done only for calculation of the criteria and for result imaging and is not used in algorithm iterations.

The accuracy of the wavefront reconstruction is characterized by $RMSE$ criteria calculated independently for phase and amplitude:

$$RMSE_\varphi = \sqrt{\frac{1}{n}||\hat{\varphi}_{est} - \varphi_{true}||_2^2},$$

$$RMSE_{ampl} = \sqrt{\frac{1}{n}||\hat{a}_{est} - a_{true}||_2^2}$$

where $n$ is a size of the image in pixels.

The noise level in the Poissonian observations is controlled in (6) by the exposure-time parameter $\chi$ and can be characterized by SNR (8) and by the mean value of photons per pixel (9).

Variations of SNR and N$_{photon}$ depend on amplitude and phase images. In our experiments with $\chi \in [0.00001, 1]$, these variations approximately take values from $-1dB$ to $60dB$ for $SNR$ and from $0.5$ to $6 \cdot 10^4$ for the mean value of photons per pixel.

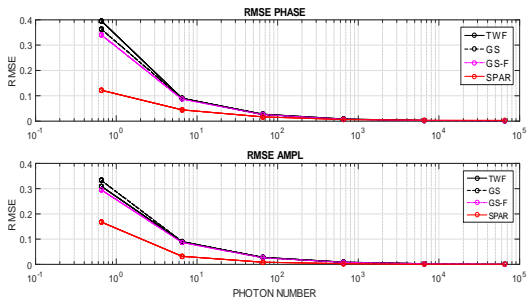

Fig. 2. Lena phase image: RMSE for phase and amplitude reconstructions versus mean number of photons per pixel. Comparison of the four algorithms: TWF, GS, GS-F and SPAR, $S = 12$.

In what follows we compare the four algorithms: TWF, GS, GS-F and SPAR. Results presented in this section can be reproduced by running the publicly available MATLAB demo-codes[1].

### A. Processing without phase unwrapping (principal phase imaging)

In experiments for principal phase imaging we show the results for two test-images: Lena ($256 \times 256$) and USAF chart ($620 \times 620$). These images are scaled to the interval $[0, \pi/2]$ to be used as an object phase. The phase unwrapping is not required for these experiments and Step 4 is omitted in the SPAR algorithm. The amplitudes

[1] http://www.cs.tut.fi/~lasip/DDT/index3.html

of the complex-valued images are invariant and equal to 1.

The results in $RMSE$ values are shown in Figs. 2 and 3, respectively, for the phase images Lena and USAF chart.

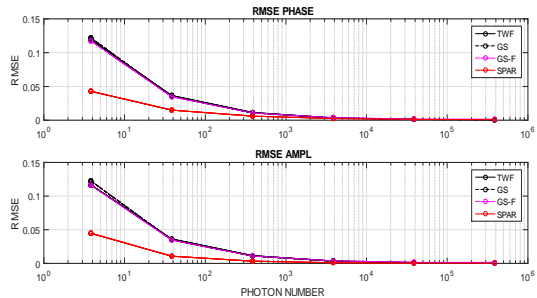

Fig. 3. USAF chart phase image: RMSE for phase and amplitude reconstructions versus mean number of photons per pixel. Comparison of the four algorithms: TWF, GS, GS-F and SPAR, $S = 12$.

The valuable advantage of SPAR for the noisier data (lower number of photons) is obvious for the mean number of photons per pixel up to 100 for Lena and up to 1000 for USAF. For larger values of photons all algorithms are about equivalent in the accuracy and demonstrate a nearly perfect reconstruction. Comparing TWF, GS and GS-F for the noisy data note that the best performance is shown by GS-F, GS shows a bit lower accuracy and the worst one in this group is the performance of TWF. Thus, SPAR is a definite favorite and other algorithm are ordered as GS-F, GS, TWF. This conclusion about the comparative accuracy of the algorithms is confirmed by our experiments for various phase images and various parameters of the experiments.

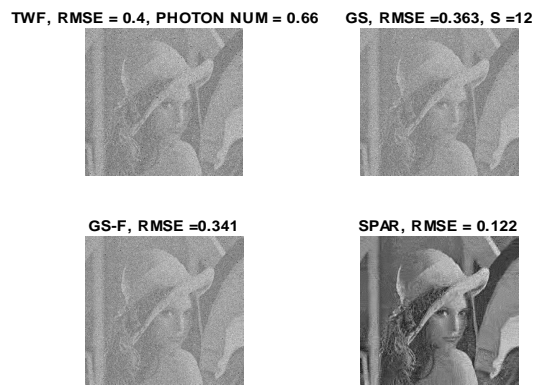

Fig. 4. Lena phase imaging for noisy data. Comparison of the TWF, GS, GS-F and SPAR algorithms, $S = 12$.

We show the reconstructed images for Lena and USAF for the most noisy cases in Figs. 4 and 5. It is clear that visually SPAR provides the best result for the Lena phase in Fig. 4. The reconstructions by the other algorithms are more or less of the equal quality. Numerically, the ordering of the algorithms corresponds to the given above conclusion from the best to the worst: GS-F, GS, TWF.

The quality of imaging for USAF in Fig. 5 overall is again with the clear advantage of SPAR, visually and numerically.

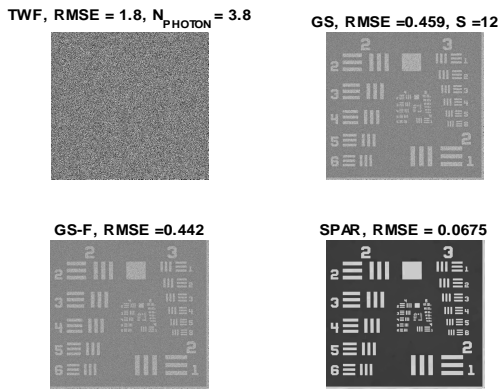

Fig. 5. USAF chart phase imaging for noisy data. Comparison of the TWF, GS, GS-F and SPAR algorithms, $S = 12$.

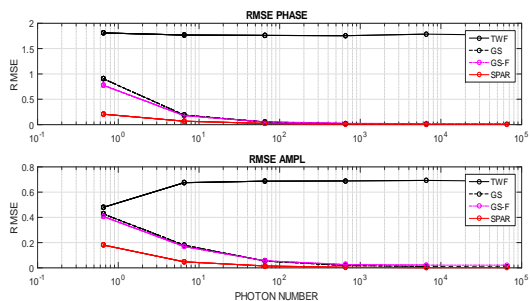

Fig. 6. Lena phase image: RMSE for phase and amplitude reconstructions versus mean number of photons. Comparison of the four algorithms: TWF, GS, GS-F and SPAR, $S = 4$.

The results in Figs. 2-5 are obtained for $S = 12$, i.e. for 12 experiments with different phase modulation masks.

For comparison we show the results obtained for much lower number of experiments $S = 4$, Figs. 6-8. The other parameters are identical to those in the above experiments. One of the important conclusions, TWF fails for this number of the experiments even for noiseless data with a very larger number of photons per pixel. An additional study shows that TWF is more sensitive to the number of the experiments as compare with other algorithms and becomes successful starting from $S \geq 10$. Figs. 6-9 confirm a strong advantage of SPAR for noisy data.

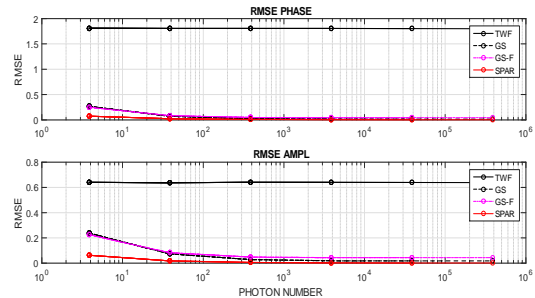

Fig. 7. USAF chart phase image: RMSE for phase and amplitude reconstructions versus the mean number of photons per pixel. Comparison of the four algorithms: TWF, GS, GS-F and SPAR, $S = 4$.

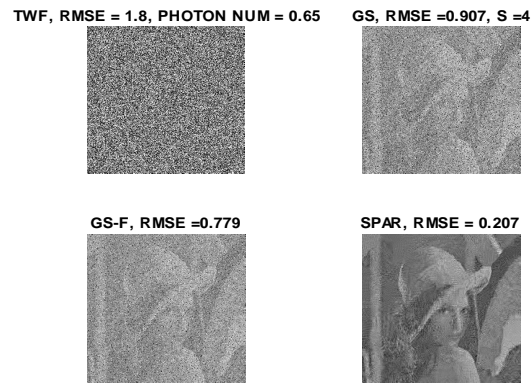

Fig. 8. Lena phase imaging for noisy data. Comparison of the TWF, GS, GS-F and SPAR algorithms, $S = 4$.

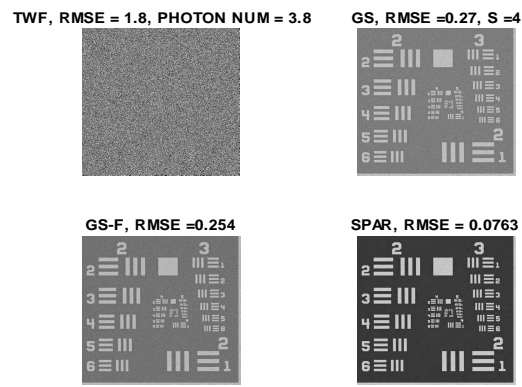

Fig. 9. USAF chart phase imaging for noisy data. Comparison of the TWF, GS, GS-F and SPAR algorithms, $S = 4$.

## B. Absolute phase imaging with phase unwrapping

In experiments for absolute phase imaging we show the results for three complex-valued data sets of size $100 \times 100$ with the invariant amplitude equal to 1 and spatially varying absolute phase: Gaussian (phase range 44 radians), truncated Gaussian (phase range 44 radians), and Shear Plane (phase range 149 radians), see Fig.10. The multiple fringes of the wrapped phases and discontinuities of the absolute phase for the truncated Gaussian and Shear plane demonstrate how complex are these tests for reconstruction, in particular, because the observations are defined by the wrapped phases. It is clear also that the sparse modeling is much simpler for the piece-wise smooth absolute phases than for the corresponding wrapped phases.

For the phase imaging we apply the four considered algorithms: TWF, GS, GS-F and SPAR. The first three algorithms give the wrapped phase reconstructions which are unwrapped by the PUMA algorithm. The SPAR algorithm is used with Step 4 (unwrapping) and, then, the sparsification and filtering are applied to the absolute phases. The absolute phase is a natural output of this algorithm.

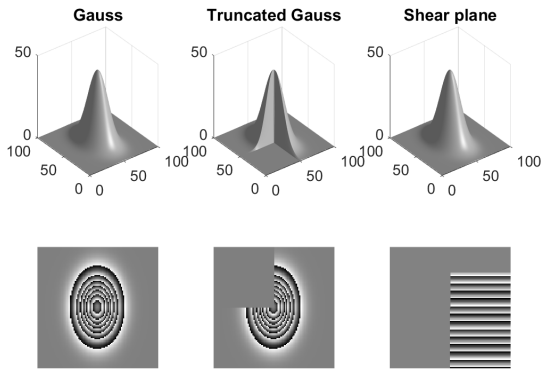

Fig. 10. True absolute phase test images: in the first row - the absolute phases, in the second row - the corresponding wrapped phases.

We do not show the RMSE curves for these experiments as qualitatively they are quite similar to those discussed above and cannot change the conclusion about the comparative accuracy of the algorithms. We demonstrate only 3D reconstructed phase surfaces obtained for the very noisy data, $N_{photon}=0.2$, because in this case the difference in the algorithm performance appears to be clear (see Figs. 11-13).

Note, that the results for the tests with for the low noise level (not shown in the paper) are more or less identical with a nearly perfect reconstruction of the phase and amplitude.

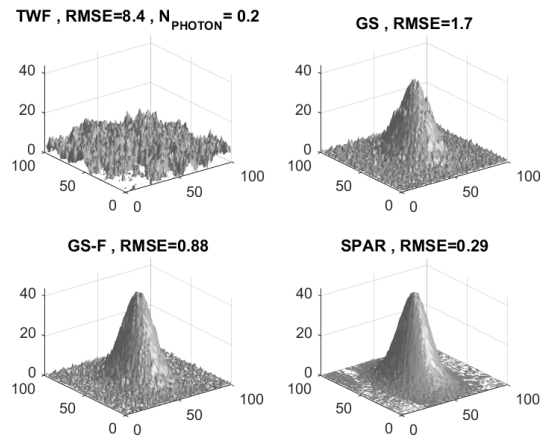

Fig. 11. The absolute (unwrapped) phase imaging for the Gaussian phase object from very noisy Poissonian observations, $S=12$.

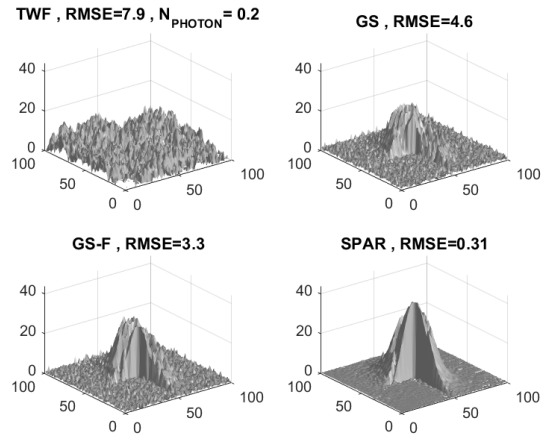

Fig. 12. The absolute (unwrapped) phase imaging for the truncated Gaussian phase object from very noisy Poissonian observations, $S=12$.

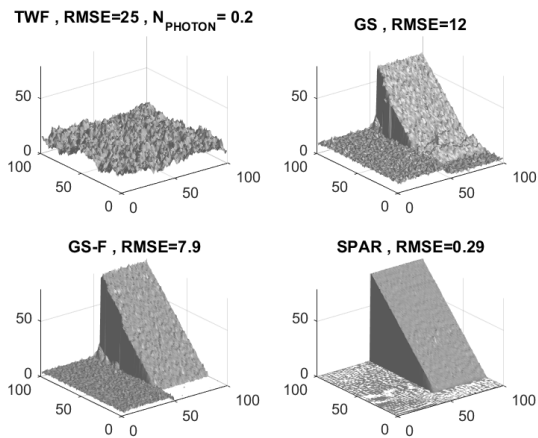

Fig. 13. The absolute (unwrapped) phase imaging for the shear plane phase object from very noisy Poissonian observations, $S=12$.

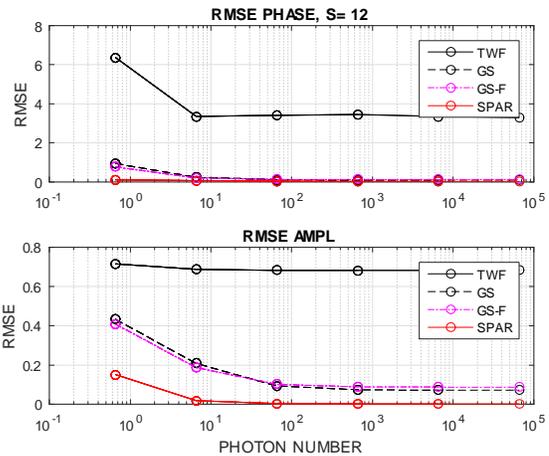

Fig. 14. Undersampled data, $p=25\%$. Lena phase image: RMSE for phase and amplitude reconstructions versus mean number of photons per pixel. Comparison of the four algorithms: TWF, GS, GS-F and SPAR, $S=12$.

## C. Undersampled data

It is a common situation in optics that only a part and even a small part of the diffraction pattern is registered by the sensor. It is a central part of the diffraction pattern corresponding to the lower frequency observations. The phase retrieval problem is to reconstruct the high-resolution image from these undersampled data.

In this scenario the criterion $\mathcal{L}_1$ in (11) takes the form

$$\mathcal{L}_1(\mathbf{u}_s, \mathbf{v}_s) = \sum_{s=1}^{S} \sum_{l\in\Omega} [|\mathbf{u}_s[l]|^2 \chi - \mathbf{z}_s[l] \log(|\mathbf{u}_s[l]|^2 \chi)]$$
$$+ \frac{1}{\gamma_1} \sum_{s=1}^{S} ||\mathbf{u}_s - \mathbf{v}_s||_2^2, \quad (34)$$

where $\Omega$ denotes a set of pixels of the diffraction pattern sampled (registered) by the sensor.

Minimization of $\mathcal{L}_1(\mathbf{u}_s, \mathbf{v}_s)$ with respect $\mathbf{u}_s$ resulting in Step 2 of the SPAR algorithm gives the solution in form

$$\hat{\mathbf{u}}_s[l] = \begin{cases} \mathbf{b}_s[l] \exp(j\cdot angle(\mathbf{v}_s[l])), & \text{if } l \in \Omega, \\ \mathbf{v}_s[l] & \text{if } l \notin \Omega, \end{cases} \quad (35)$$

where $\mathbf{b}_s[l]$ is defined by (20) for the Poissonian data.

Here, the vectors $\mathbf{u}$ and $\mathbf{v}$ have the same dimension, $\mathbf{u}, \mathbf{v} \in \mathbb{C}^n$, and only $\mathbf{u}_s[l]$ with $l \in \Omega$ are subjects of the amplitude update due to the given observations while others are equal to $\mathbf{v}_s[l]$, $l \notin \Omega$.

This extrapolation of $\hat{\mathbf{u}}_s$ outside of the sensor area is used also in our experiments with the GS and GS-F algorithms.

Let $p$ denote a percentage of the diffraction pattern pixels used for processing. Simulation experiments show that the modified in this way GS, GS-F and SPAR algorithms works very well provided $p \geq 25\%$ for quite noisy data, while the standard versions of these algorithms, where $\hat{\mathbf{u}}_s[l]$ is zeroed for $l \notin \Omega$, actually fail even for noiseless data. It proves that the modification of the algorithms defined by (35) is of importance.

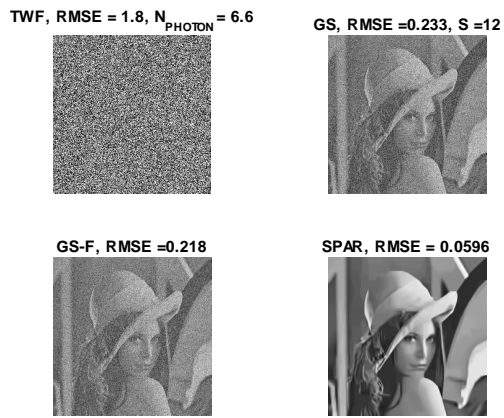

Fig. 15. Undersampled data, $p=25\%$. Lena phase imaging for noisy data. Comparison of the TWF, GS, GS-F and SPAR algorithms, $S=12$.

Figs. 14-15 produced for $p=25\%$ and the Lena phase image show the RMSE curves and image reconstruction from noisy data. One may note that TWF fails in this scenario. However, this comparison is not completely fair as the extrapolation rule (35) is not used in this algorithm. Instead of this extrapolation it works with zeroing of $75\%$ of the diffraction pattern area outside of its central part.

In Fig.16 we show the reconstructions of the USAF phase test image from noisy data for $p$ as small as $p = 16\%$. The GS and GS-F algorithms are nearly failed in this case. Contrary to it, the SPAR algorithm is quite successful.

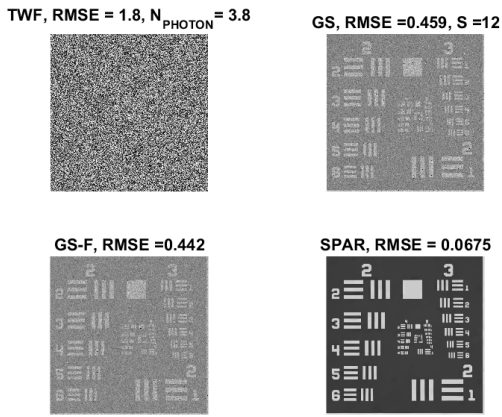

Fig. 16. Undersampled data, $p = 16\%$. USAF chart phase imaging for noisy data, noisy data. Comparison of the TWF, GS, GS-F and SPAR algorithms, $S = 12$.

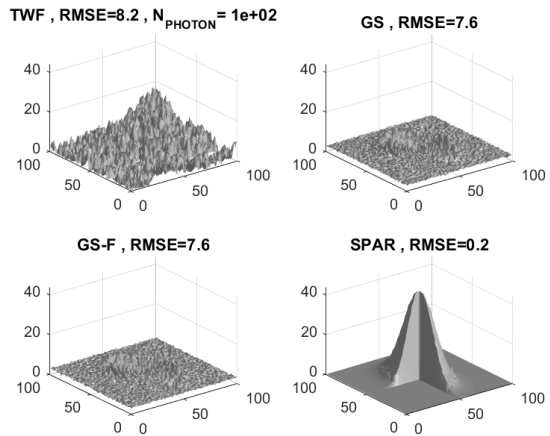

Fig. 18. Undersampled data, $p = 16\%$. Truncated Gaussian absolute phase imaging. Comparison of the TWF, GS, GS-F and SPAR algorithms, $S = 12$.

Let us apply the algorithms developed for the undersampled data to the absolute phase test-images from Subsection IV-B. In what follows, the results are shown for the data with $N_{photon} = 100$ which are quite noisy for these tests-images and $p$ as small as $p = 16\%$. The reconstructions produced by the four compared algorithms are shown in Figs.14-19. We may note that all algorithms except SPAR are failed and failed completely being not able to reveal the 3D phase shapes of the objects. The quality of the SPAR reconstructions visually and numerically is quite acceptable.

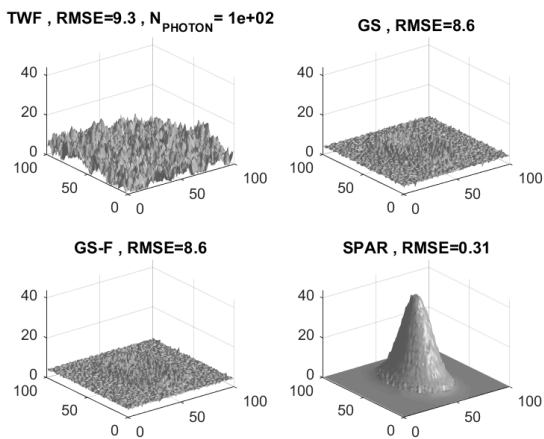

Fig. 17. Undersampled data, $p = 16\%$. Gaussian absolute phase imaging. Comparison of the TWF, GS, GS-F and SPAR algorithms, $S = 12$.

It is important to emphasize that in the above experiments with the principal phase objects the advantage of the SPAR algorithm concerns only the accuracy of reconstructions and lower values of RMSE, while in the considered experiments with the absolute phase test-image the advantage of SPAR is of a principal nature as the counterpart algorithms failed to solve the problem while SPAR is quite successful.

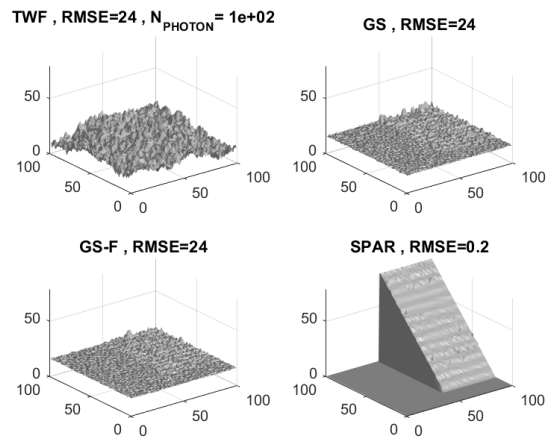

Fig. 19. Undersampled data, $p = 16\%$. Shear plane absolute phase imaging. Comparison of the TWF, GS, GS-F and SPAR algorithms, $S = 12$.

Surprisingly, SPAR being quite successful for $p = 16\%$ and $N_{photons} = 100$ was not able to improve the accuracy of reconstruction for the much lower level noise and even for noiseless data. We assume that in this scenario the main problem is the accurate object modeling with interpolation and extrapolation of variables while the denoising is a secondary problem.

## D. Parameters of the SPAR algorithm

The performance of the SPAR algorithm essentially depends on its parameters. In our experiments the parameters are fixed for all tests. The image patches in $BM3D$ are square $8 \times 8$. The group size is limited by $25$ patches. The step size between the neighboring patches is equal to $3$. The transforms DCT (for patches) and Haar (for the group length) are used for 3D group data processing in $BM3D$. In the shown results as an initial guess for the iterative GS and SPAR algorithm we use an image with the invariant amplitude equal to $1$ and the zero mean Gaussian i.i.d. random phase with the standard deviation $0.1\pi$. The number of the iterations is fixed to $50$.

The parameters defining the iterations of the algorithm are as follows: $\gamma_1 = 1/\chi$; $th_a^0 = 1.4$; $th_\varphi^0 = 1.4$. For the experiments with the undersampled data we use the larger threshold values $th_a^0 = 5.6$; $th_\varphi^0 = 5.6$.

The complexity of the algorithm's iterations are defined by the built-in BM3D filters generating data adaptive synthesis and analysis matrices varying in iterations. The theoretical analysis of the complexity of BM3D can be seen in [35], in particular, the computational time is proportional to the image size provided that the other parameters of the algorithm are fixed.

For our experiments we use MATLAB R2014a and the computer with the processor Intel(R) Core(TM) i7-4800MQ@ 2.7 GHz. The complexity of the algorithm is characterized by the time required for processing. For $50$ iterations and $256 \times 256$ images this time is as follows: TWF$\simeq$ 8.7 sec.; GS$\simeq$ 3.7 sec.; GS-F$\simeq$ 4.3 sec., SPAR$\simeq$ 65 sec. (no unwrapping); SPAR$\simeq$ 90 sec. (with unwrapping).

## V. CONCLUSION

The phase retrieval from intensity observations is considered. This paper introduces a variational approach to object phase and amplitude reconstruction from noisy Poissonian intensity observations. The maximum likelihood criterion used in the developed multiobjective optimization (Nash equilibrium technique) defines the intention to reach statistically optimal estimates. The sparsity is one of the key elements of the developed SPAR algorithm used for modeling of spatially varying amplitude and phase. The phase retrieval is an ill-possed inverse problem where the observation noise is amplified and transferred to phase and amplitude as variables of optimization. The sparse modeling enables a regularization of this inverse problem and efficient suppression of these random errors by BM3D filtering of phase and amplitude. The efficiency of the SPAR algorithm is demonstrated by simulation experiments for the coded diffraction pattern scenario. The comparison is produced versus the TWF, GS-F and GS algorithms. For the noisy observations the SPAR algorithm demonstrates a serious advantage. For the low noise level the accuracy of the SPAR algorithm as well as its simplified versions the GS and GS-F algorithms is nearly identical to the accuracy of the TWF algorithm. The GS and GS-F algorithms are essentially faster than TWF while SPAR computationally much more demanding is slower than TWF, GS and GS-F. SPAR derived as an optimal algorithm for the Poissonian observations due to the strong BM3D filters is robust with respect to different kind of additional errors such as an additive Gaussian noise in observations and quantization errors. This robustness was tested in our study which is not included in this paper.

## VI. ACKNOWLEDGEMENT

This work was supported by Academy of Finland, project no. 287150, 2015-2019.